\documentclass[preprint,authoryear,10pt]{elsarticle}
\newtheorem{theorem}{Theorem}

\newtheorem{corollary}{Corollary}[theorem]
\newtheorem{example}{Example}
\newtheorem{remark}{Remark}
\newtheorem{lemma}{Lemma}
\newproof{pf}{Proof}
\newtheorem{assum}{Assumption}

\usepackage{amssymb}

\journal{Linear Algebra and its Applications}

\begin{document}

\begin{frontmatter}



\title{Yet another approach to the Algebraic Riccati Inequality}


\author[rvt]{A. Sanand Amita Dilip\corref{cor1}}
\ead{sanand@ee.iitkgp.ac.in}
\author[focal]{Harish K. Pillai\corref{cor2}}
\ead{hp@ee.iitb.ac.in}

\address[rvt]{Department of Electrical Engineering, IIT Kharagpur, Kharagpur, India}
\address[focal]{Department of Electrical Engineering, IIT Bombay, Powai, India}

\begin{abstract}
 We give a rank characterization of the solution set of algebraic Riccati inequality (ARI) for both controllable and
 uncontrollable systems. Assuming an existence of a solution 
 of the corresponding algebraic Riccati equation (ARE), we characterize the boundedness/unboundedness properties of solutions of ARI for controllable/uncontrollable systems 
 without any assumption on sign controllability. As a consequence of our observations, we obtain Willems' result 
 $K_{min}\leq K\leq K_{max}$  for an ARI in the case of controllable systems and explore some structure on the extremal solutions. 
 We also consider the curious case of 
 uncontrollable purely imaginary eigenvalues and the behavior of the solution set of ARI. 
 In particular, we show that a system is controllable if and only if the set of solutions of an ARI is bounded. In addition, we study the effect of the position 
 of eigenvalues of the system matrix in the complex plane on the behavior of the solution set of ARIs. 
 Furthermore, we obtain a rank parametrization for solutions of ARI for controllable systems. 
\end{abstract}


\end{frontmatter}


\section{Introduction}
Algebraic Riccati inequality (ARI) arises in $H_{\infty}$ control (\cite{schererV}, \cite{schererIV}, \cite{schererI}) and also in the formulation of storage functions
 for dissipative systems (\cite{willems}).  
Study of symmetric solutions of ARI has appeared in \cite{schererII}, 
\cite{schererIII}, \cite{fau}, \cite{lind1}, \cite{lind2}, \cite{ferrIV}, \cite{pav} and some of the references therein. 
 Our focus is on the symmetric ARIs of the form 
\begin{eqnarray}
  -A^TK - KA -Q + KBB^TK\leq 0\nonumber
 \end{eqnarray}
 where
 $A,B,Q$ are real matrices having dimensions $n\times n,n\times m$ and $n\times n$ respectively with 
 $Q$ and $K$ being  symmetric. (Note that in the characterization of storage functions, the ARI takes the
form $Q - A^TK-KA- KBB^TK \geq 0$ which one obtains from the ARI: $-A^TK - KA -Q + KBB^TK \leq 0$ by replacing $K$ by $-K$.)
The solution set of the ARI (which is a spectrahedron) characterizes the set of all possible storage functions for dissipative systems (\cite{willems}). 
 Characterization of solutions of ARI for sign-controllable 
  $(A,B)$ pairs is well studied in the literature (\cite{faib1}, \cite{schererII}). 
  For necessary and sufficient conditions for existence of a solution for 
  non-strict ARIs, we refer the reader to \cite{schererIII}. 

We primarily consider homogeneous ARIs in this paper. A non-homogeneous ARI can be converted into a homogeneous one by fixing a solution of the 
correspond ARE \cite{willems}, \cite{schererII}, \cite{schererI}, \cite{MR945585}. 
 For a detailed study of continuous time homogeneous algebraic Riccati equations and their applications, we refer the reader to \cite{ferrI}, \cite{ferrII}, 
\cite{ferrIII}, \cite{ferrIV}, \cite{picc} and the references therein.
For non-homogeneous AREs, we refer the reader to \cite{willems}, \cite{lanc}, \cite{MR1133203}, 
\cite{mehr2}, \cite{mehr3}, \cite{vandooren}, \cite{wim2}, \cite{wim3}.
 Reader may refer \cite{bittan}, \cite{MR2896454}, \cite{MR1997753}, 
 \cite{MR1133203} for ARE and their applications. 

We gave a rank characterization of solutions of symmetric AREs in \cite{are} and non-symmetric AREs is \cite{nare}. Then we extended our approach to 
 discrete time AREs in \cite{dare}. In here, we extend our approach in these previous works to obtain a rank characterization of solutions of continuous time homogeneous ARIs. 
 This in turn provides a unified linear algebraic perspective to AREs (symmetric and non-symmetric, continuous and discrete) and ARIs.
 

  Note that sign-controllability of $(A,B)$ pair rules out purely imaginary uncontrollable 
  eigenvalues of the feedback matrices. With the following assumption (Assumption \ref{assum1}), we allow feedback matrices to have purely imaginary eigenvalues.   

\begin{assum}\label{assum1}
 We assume that a solution to $-A^TK-KA-Q+KBB^TK=0$ always exists.
\end{assum}

We 
 fix an arbitrary solution $K_0$ of the ARE $-A^TK-KA-Q+KBB^TK=0$.
Let $K = K_0 + X$ where $X$ can be thought of as a perturbation from $K_0$. A similar construction of fixing a solution is used in \cite{willems}, \cite{schererII}, \cite{schererI}, \cite{MR945585}. 
We can then re-write $-A^TK-KA-Q+KBB^TK$ as
 \begin{eqnarray}
   &=& -A^T(K_0 + X) - (K_0 + X)A -Q + (K_0 + X)BB^T(K_0 + X)\nonumber\\
 &=& -A^TK_0 - K_0A -Q + K_0BB^TK_0 -A^TX - XA + K_0BB^TX\nonumber\\
&& + XBB^TK_0 + XBB^TX
\nonumber\\
 &=&  - (A - BB^TK_0)^TX - X(A - BB^TK_0) + XBB^TX.\label{eq1}\\
 && (\mbox{Since } -A^TK_0 - K_0A -Q + K_0BB^TK_0 = 0.) \nonumber
\end{eqnarray}
Let $A_0 = A-BB^TK_0$. 
We denote $- A_0^TX - XA_0 + XBB^TX$ by the notation Ric$(X): = -A_0^TX-XA_0 + XBB^TX$.
Note that 
we are interested in real, symmetric solutions of Ric$(X)\leq0$.
  We characterize the solutions of Ric$(X)\leq0$ in terms of eigenspaces/invariant subspaces of the newly constructed matrix $A_0^T$ where we obtain a 
  criterion for the boundedness 
  of the solution set and give a rank parametrization of the solution set. 
 
 
 {\bf Paper Organization}: The paper is organized as follows. In Section $2$, we build some preliminaries to be used in the paper followed by Section 
 $3$; where we consider degenerate cases of zero eigenvalue and purely imaginary eigenvalues of the feedback matrix $A_0$ and characterize the 
 solution set of ARI for these specific cases.  
 Our main results are stated in sections \ref{ari} and \ref{para}. In Section \ref{ari}, we study 
 boundedness/unboundedness properties of solutions of ARIs for a completely general case when $(A,B)$ is controllable/uncontrollable and 
 Section \ref{para} involves a rank parametrization 
 of solutions of ARIs. 
 \subsection{Notation} We briefly mention the terminology followed in the rest of the paper. All the matrices considered here are real. They may have real or 
 complex eigenvalues which may be distinct or repeated. By $D_J$ we denote the upper triangular Jordan canonical form associated with $A_0^T$ where 
 $A_0=A-BB^TK_0$ (Equation \ref{eq1}). 
 By an order of an equation we mean the size of the matrices involved in 
 that equation. We list some of the frequently appearing equations in the paper as follows: 
 \begin{eqnarray}
  &&\mbox{Ric}(X)=-A_0^{T}X-XA+XBB^TX.\nonumber\\
  &&\mbox{Simplified ARE}: -D_J\mathcal{L}-\mathcal{L}D_J^T+\mathcal{L}M\mathcal{L}=0.\nonumber\\
  &&\mbox{Reduced order simplified ARE}:-D_{J_k}\mathcal{L}_k-\mathcal{L}_kD_{J_k}^T+\mathcal{L}_kM_k\mathcal{L}_k=0.\nonumber\\
  &&\mbox{Simplified ARI}: -D_J\mathcal{L}-\mathcal{L}D_J^T+\mathcal{L}M\mathcal{L}\le0.\nonumber\\
  &&\mbox{Reduced order simplified ARI}: -D_{J_k}\mathcal{L}_k-\mathcal{L}_kD_{J_k}^T+\mathcal{L}_kM_k\mathcal{L}_k\leq0.\nonumber
 \end{eqnarray}

\section{Preliminaries}\label{bldblk}
Given a pair of matrices $(A_0,B)$, using a change of basis, it is possible to write $A_0$ and $B$ in the following form (\cite{won}, \cite{kai})
\begin{eqnarray}
 A_0=\left[\begin{array}{cc}A_0^{11}&A_0^{12}\\0&A_0^{22}\end{array}\right], B= \left[\begin{array}{c}B_{1}\\0\end{array}\right]\nonumber
\end{eqnarray}
 where $(A_0^{11},B_1)$ form a controllable pair. Eigenvalues of $A_0^{11}$ are controllable while eigenvalues of $A_0^{22}$ are
 uncontrollable. Left eigenvectors of $A_0$ (which are right eigenvectors of $A_0^T$) corresponding to the uncontrollable eigenvalues are of the form
 $v^T=\left[\begin{array}{cc}0&u^T\end{array}\right]$ and $v^TB=0$. 
 All right
 eigenvectors of $A_0^T$ that belong to the kernel of $B^T$, are called uncontrollable modes. If $v$ is an eigenvector of $A_0^T$ associated with a
 controllable eigenvalue, then $B^Tv\neq0$ and these eigenvectors correspond to controllable modes. We use these facts to characterize the solution 
 set of ARIs for controllable and uncontrollable systems. 

 Consider the ARI Ric$(X)=-A_0^TX-XA_0 + XBB^TX \le0$. 
 Clearly, $X =0$ is a trivial solution of Ric$(X) \le0$. We now look for a nonzero $X$ that satisfies Ric$(X) \leq 0$. 
We begin with a simplest case i.e. matrices $X$ that have rank one.  Since $X$ is symmetric, let $X = \alpha vv^T$ where $\alpha \in \mathbb{R}$
and $v \in \mathbb{R}^{n}$ with $||v|| = 1$.  (Note that by $||.||$, we mean the two norm of a vector.) 
\begin{theorem}\label{rk1lem}
  Let $X=\alpha vv^T$, such that $||v||=1$. 
  \begin{enumerate}
   \item If $v$ is an eigenvector of $A_0^T$ then, the rank of Ric$(X)$ is at most one and Ric$(X)$ is semi-definite.
   \item If $v$ is not an eigenvector of $A_0^T$ then, Ric$(X)$ is indefinite.
  \end{enumerate}

\end{theorem}
\begin{pf}
 Refer Theorem $1$ of \cite{are}. \qed
\end{pf}
 The following theorem is a rank two analogue of Theorem \ref{rk1lem}. 
 It characterizes all the allowable rank two perturbations from the fixed solution such that the ARI Ric$(X)\leq 0$ gets satisfied.
\begin{theorem}\label{rk2per}
 If  $X=L\mathcal{L}L^T$ (where $L$ is $n\times2$ and $\mathcal{L}$ is $2\times2$) such that two columns of $L$ are linearly independent,
 then Ric$(X)$ is semi-definite only if the column span of $L$ is a $A_0^T-$invariant subspace.
\end{theorem}
\begin{pf}
Refer Theorem $4$ of \cite{are}.\qed
\end{pf}

  This theorem holds for general rank $k$ perturbations also and the proof can be given along similar lines. Thus, we assume that $X$ is of the form 
  $X=L\mathcal{L}L^T$ where the column span of $L$ forms an $A_0^T-$invariant subspace. 


\begin{lemma}\label{smplfd}
 Let $X = L\mathcal{L}L^T$ where $A_0^TL=LD_J$. Then, $X$ is a solution of Ric$(X)\le0$ if and only if 
 $\mathcal{L}$ satisfies the simplified ARI 
 \begin{eqnarray}
  -D_J\mathcal{L}-\mathcal{L}D_J^T+\mathcal{L}M\mathcal{L}\le0.\label{smfdari}
 \end{eqnarray}

\end{lemma}
\begin{pf}
 If $\mathcal{L}$ satisfies the simplified ARI $-D_J\mathcal{L}-\mathcal{L}D_J^T+\mathcal{L}M\mathcal{L}\le0$, then clearly, Ric$(X)\le0$. 
 Conversely, if Ric$(X)\le0$, then the column span of $X$ must be $A_0^T-$invariant. Since $X$ is symmetric, $X = L\mathcal{L}L^T$ with the 
 column span of $L$ forming an $A_0^T$ invariant subspace. Thus, Ric$(X)\le0$ simplifies to 
 $-D_J\mathcal{L}-\mathcal{L}D_J^T+\mathcal{L}M\mathcal{L}\le0$ and the lemma follows. \qed
\end{pf}
\begin{remark}
Let $L \in \mathbb{R}^{n\times n}$ be such that $A_0^TL=LD_J$ 
where $D_J$ is the Jordan form (upper triangular) associated with $A_0^T$. 
Using $X = L\mathcal{L}L^T$ (where $\mathcal{L}$ is $n\times n$ matrix),
Ric$(X)$ is reduced to an expression
$L(-D_J\mathcal{L}-\mathcal{L}D_J^T+\mathcal{L}M\mathcal{L})L^T$ where $M=L^TBB^TL$. 
Therefore, solving $-D_J\mathcal{L}-\mathcal{L}D_J^T+\mathcal{L}M\mathcal{L}\le0$ is equivalent to solving Ric$(X)\le0$.
\end{remark}
\begin{remark}
We observe that any rank $k$ solution of Ric$(X) \le 0$ can be written as $X = L_k\mathcal{L}_kL_k^T$ where $\mathcal{L}_k$ is a 
symmetric $k \times k$ matrix  and the column span of $L_k$ ($L_k \in \mathbb{R}^{n\times k}$) forms an $A_0^T-$invariant subspace. 
Using $X=L_k\mathcal{L}L_k^T$, the ARI Ric$(X) \le 0$ reduces to the form $-D_{J_k}\mathcal{L}_k-\mathcal{L}_kD_{J_k}^T+\mathcal{L}_kM_k\mathcal{L}_k\leq0$
(where $M_k=L_k^TBB^TL_k$ and $D_{J_k}$ is a $k \times k$ Jordan sub-block). Thus, by the above lemma, $\mathcal{L}_k$ must satisfy the simplified 
ARI above involving $k\times k$ matrices. We refer to the ARI $-D_{J_k}\mathcal{L}_k-\mathcal{L}_kD_{J_k}^T+\mathcal{L}_kM_k\mathcal{L}_k\leq0$ as 
the reduced order simplified ARI.
\end{remark}

The above lemma and remarks allow us to simplify the structure of ARIs using invariant subspaces of $A_0^T$. When $A_0$ has zero or purely imaginary eigenvalues, the behavior of the solution 
set of ARI is slightly different for controllable systems. In the following section, we consider these cases for both controllable and uncontrollable systems. 

\section{Zero and purely imaginary eigenvalues of $A_0$}\label{genrk}
We observe that the invariant subspaces associated with zero and purely imaginary eigenvalues of $A_0$ form degenerate cases for controllable systems. 
We also consider the interesting case of uncontrollable purely imaginary eigenvalues (Corollary $\ref{uncpuim}$) which has not been studied in the literature to the best of our knowledge.
 
 We give below a lemma which will be used repeatedly in the subsequent results. 
\begin{lemma}\label{indeflem}
Suppose $P\geq0$, $Q$ is indefinite and $P+Q = R$. Then, $R$ can not be negative semidefinite.
\end{lemma}
\begin{pf}
 Suppose $R\leq0$ where $P+Q = R$.
It is clear that $P-R = -Q$ but $P$ and $-R$ are positive semidefinite and $-Q$ is indefinite which is a contradiction. Thus, $R$ cannot be $\leq0$.
\end{pf}
\subsection{Zero eigenvalues of $A_0$}
Consider the first degenerate case for controllable systems where $A_0$ has a zero eigenvalue.
\begin{theorem}\label{zerocont}
 Suppose $(A_0,B)$ is controllable and $A_0$ has only one eigenvalue which is equal to zero. Let 
 $X=L_j\mathcal{L}_j L_j^T$ where columns of $L_j$ form generalized eigenvectors corresponding to zero eigenvalue of 
 $A_0^T$.  Then, there are no non-trivial solutions of Ric$(X)\le0$. 
\end{theorem}
\begin{pf}
 Suppose that the columns of $L_j$ form linearly independent eigenvectors of $A_0^T$ corresponding to the zero eigenvalue. Hence, $A_0^TL_j=0$. Therefore, for $X=L_j\mathcal{L}_jL_j^T$, Ric$(X) = L_j\mathcal{L}_jM_j\mathcal{L}_jL_j^T \geq0$ 
 where $M_j = L_j^TBB^TL_j$. Now since zero is a controllable eigenvalue, $L_j^TB \neq 0$. Therefore, $M_j\neq 0$ 
 and Ric$(X)\leq0$ has no non-zero solutions.\\
 \indent Now suppose columns of $L_j$ denote generalized eigenvectors of $A_0^T$. For simplicity, assume that there is only one Jordan block 
 for the zero eigenvalue. $A_0^TL_j=L_jD_{J_j}$ where all diagonal entries of $D_{J_j}$ are zero. For $X=L_j\mathcal{L}_jL_j^T$, Ric$(X) = 
 L_j(-D_{J_j}\mathcal{L}_j-\mathcal{L}_jD_{J_j}^T+\mathcal{L}_jM_j\mathcal{L}_j)L_j^T$. Let $\mathcal{L}_j = 
 \left[\begin{array}{cccc}l_1&l_2&\cdots&l_j\end{array}\right]$. Therefore, 
 \begin{eqnarray}
  -D_{J_j}\mathcal{L}_j-\mathcal{L}_jD_{J_j}^T = 
 -\left[\begin{array}{ccccc}l_2&l_3&\cdots&l_j&0\end{array}\right]-\left[\begin{array}{c}l_2^T\\l_3^T\\.\\.\\l_j^T\\0\end{array}\right]=C.\nonumber
 \end{eqnarray}

  Consider the determinant of $2\times2$ principal sub-matrix obtained from the first and the $j-$th rows and columns respectively of this matrix $C$. This determinant is negative which 
 implies $-D_{J_j}\mathcal{L}_j-\mathcal{L}_jD_{J_j}^T$ is indefinite. But $\mathcal{L}_jM_j\mathcal{L}_j$ is positive semidefinite. Hence, 
 Ric$(X)$ cannot be negative semidefinite unless $X=0$ (Lemma $\ref{indeflem}$). Hence, the ARI is not satisfied for any non-zero $X$. 
 
 The general case when there 
 are more than one Jordan blocks also follows similarly. 
Suppose there exist a full rank solution of $-D_{J_j}\mathcal{L}_j-\mathcal{L}_jD_{J_j}^T+\mathcal{L}_jM_j\mathcal{L}_j\le0$ and let $Y$ be its inverse. 
Hence, $-YD_{J_j}-D_{J_j}^TY+M_j\leq0$. We restrict ourselves to the first Jordan block $D_{J_j}^1$ of $D_{J_j}$. We partition $Y$ and $M_j$ accordingly. 
Thus, we have $-Y^1D_{J_j}^1-(D_{J_j}^1)^TY^1+M_j^{1} \leq 0$. Since $-Y^1D_{J_j}^1-(D_{J_j}^1)^TY^1$ is indefinite and $M_j^1\geq0$, the inequality 
  $-Y^1D_{J_j}^1-(D_{J_j}^1)^TY^1+M_j^{1} \leq 0$ is possible only when $M_j^1=0$ and $Y^1=0$ (Lemma $\ref{indeflem}$). 
  But $M_j^1=0$ contradicts the controllability. Therefore, the full rank solution of $-D_{J_j}\mathcal{L}_j-\mathcal{L}_jD_{J_j}^T+\mathcal{L}_jM_j\mathcal{L}_j\le0$ can not exist. 
  Similarly, we can show that the low rank solutions do not exist by restricting to the lower dimensional invariant subspaces corresponding to the zero eigenvalue. 
  To be more precise, we may take $k-$generalized eigenvectors of $A_0^T$ forming $A_0^T-$invariant subspace as columns of $L_j$ and use 
  $X=L_j\mathcal{L}_kL_j^T$, (where $\mathcal{L}_k$ is a $k\times k$ symmetric matrix) to 
  transform  Ric$(X)\le 0$ to a $k\times k$ reduced order simplified ARI. Now using similar arguments used for the full rank case, one can show that a rank $k$ solution 
  of the ARI Ric$(X)\le 0$ does not exist. 
\qed
\end{pf}
\begin{corollary}\label{zerouncont}
 If $A_0$ has a zero eigenvalue which is uncontrollable, then $X=L_{uc}\mathcal{L}L_{uc}^T$ is a solution of Ric$(X)\leq0$ for any symmetric
 $\mathcal{L}$, where columns of $L_{uc}$ are eigenvectors of $A_0^T$ associated with the uncontrollable modes of the zero eigenvalue.
\end{corollary}
\begin{pf}
 Note that $A_0^TL_{uc} = L_{uc}^TA_0=0$ and $L_{uc}^TBB^TL_{uc}=0$ which proves the corollary.\qed
\end{pf}

\subsection{Purely imaginary eigenvalues of $A_0$}
 Now we consider the case when $A_0$ has purely imaginary eigenvalues $\pm i\mu$. 

\begin{theorem}\label{pureimag}
 Suppose $(A_0,B)$ is controllable. 
 If $X=L\mathcal{L}L^T$ where columns of $L$ form the two dimensional $A_0^T$-invariant subspace associated with a complex conjugate pair of
 purely imaginary eigenvalues $\pm i\mu$ of $A_0^T$, then both Ric$(X)\leq0$ and Ric$(X)=0$ are not satisfied for any nonzero $X$ of the given form.
\end{theorem}
\begin{pf}
We may assume that columns $v_1,v_2$ of $L$ are such that
$A_0^T\left[\begin{array}{cc}v_1&v_2\end{array}\right]= \left[\begin{array}{cc}v_1&v_2\end{array}\right]
 \left[\begin{array}{cc}0&\mu\\-\mu&0\end{array}\right]$. Let $\mathcal{L}=\left[\begin{array}{cc}a&b\\b & c\end{array}\right]$ and 
 $D=\left[\begin{array}{cc}0&\mu\\-\mu&0\end{array}\right]$.
\begin{eqnarray}
  \mbox{Ric}(X) =  L(-D\mathcal{L}-\mathcal{L}D^T+\mathcal{L}M\mathcal{L})L^T\nonumber\\
  =L\{\mu\left[\begin{array}{cc}-2b&a-c\\a-c
&2b\end{array}\right]+\left[\begin{array}{cc}a&b\\b & c\end{array}\right]
 M\left[\begin{array}{cc}a&b\\b & c\end{array}\right]\}L^T\label{eqim1}
\end{eqnarray}
\indent Note that last term is positive semidefinite since $M\geq0$. Consider \\$\mu\left[\begin{array}{cc}-2b&a-c\\a-c &2b\end{array}\right]$,
$\mu\geq0$. The determinant of this matrix is $-4b^2-(a-c)^2\leq0$. 
Suppose $-4b^2-(a-c)^2=0$. This implies that $a=c$ and $b=0$. From Equation $(\ref{eqim1})$, Ric$(X)\le 0$ if and only if $M=0$ which is not possible 
due to controllability of the $(A_0,B)$ pair. Therefore, assume that $-4b^2-(a-c)^2<0$ which implies that the matrix 
$\mu\left[\begin{array}{cc}-2b&a-c\\a-c &2b\end{array}\right]$ is indefinite. Thus, Ric$(X)$ involves the sum of a positive semidefinite
matrix and an indefinite matrix -- this sum cannot be negative semidefinite (by Lemma $\ref{indeflem}$). 
Hence, both ARE and ARI mentioned in the theorem have only the zero solution.
\qed
\end{pf}
\begin{corollary}\label{uncpuim}
 Suppose $(A_0,B)$ has a pair of purely imaginary eigenvalues which are uncontrollable. 
 If $X=L\mathcal{L}L^T$ where columns of $L$ form the two dimensional $A_0^T$-invariant subspace associated with a complex conjugate pair of
 purely imaginary eigenvalues $\pm i\mu$ of $A_0^T$, then Ric$(X)=0$ has infinitely many solutions associated with the invariant 
 subspace and the set of solutions becomes unbounded.
\end{corollary}
\begin{pf}
 Note that since the purely imaginary eigenvalues are uncontrollable, $M=L^TBB^TL=0$. Therefore, from Equation $(\ref{eqim1})$, Ric$(X)\le 0$ whenever 
 $a=c$ and $b=0$. In other words, any $X =L\left[\begin{array}{cc}a&0\\0&a\end{array}\right]L^T$ for any $a \in \mathbb{R}$ gives a solution of Ric$(X)\le0$. 
 Thus, it follows that the solution set of Ric$(X)\le0$ becomes unbounded.\qed
\end{pf}
\begin{theorem}\label{pureimagjord}
 Suppose $(A_0,B)$ is controllable. If $A_0^T$ has purely imaginary eigenvalues with trivial or non-trivial Jordan 
structure and geometric multiplicities greater than one, then 
the invariant subspace corresponding to these purely imaginary eigenvalues of $A_0^T$ do not 
contribute to solutions of Ric$(X)\leq0$.
\end{theorem}
\begin{pf}
Without loss of generality, assume that $A_0$ has a pair of purely imaginary eigenvalues repeated $\frac{n}{2}$ times where $n$ is even. 
 Suppose columns of $L$ form an invariant subspace corresponding to these purely imaginary eigenvalues. Thus, 
Ric$(X)\le0$ is reduced to $-D_J\mathcal{L}-\mathcal{L}D_J^T+\mathcal{L}M\mathcal{L}\le0$ where $D_J$ has only purely imaginary eigenvalues 
(Lemma \ref{smplfd}).  
Assuming that a full rank solution $\mathcal{L}$ exists, let $Y$ be its inverse. Therefore, pre and post multiplying $-D_J\mathcal{L}-\mathcal{L}D_J^T+\mathcal{L}M\mathcal{L}\le0$ by 
$Y$, we get $-YD_J-D_J^TY+M\le0$. Restricting to $(1,1)$ block in the expression $-YD_J-D_J^TY+M\le0$ 
determined by the position of the first purely imaginary block $D_J^1$ in Jordan structure of $D_J$, we obtain  
$-Y^1D_J^1-(D_J^1)^TY^1+M^1\le0$. This is possible only when $Y^1=0$ and $M^1=0$ (Lemma $\ref{indeflem}$). 
But $M^1\neq0$ due to controllability. Thus, $Y$ can not exist. 
One can similarly show the non-existence of 
 low rank solutions by restricting to the lower dimensional invariant subspaces corresponding to the purely imaginary eigenvalues. 
\qed 
\end{pf}

Thus, for controllable systems, invariant subspaces associated with the zero eigenvalue or purely imaginary eigenvalues do not correspond to any non-trivial 
solutions of Ric$(X)\le0$. Whereas, for uncontrollable systems where there are zero or purely imaginary uncontrollable eigenvalues, there are infinitely many 
solutions associated with the corresponding invariant subspace and the solution set becomes unbounded. These observations will be used to characterize boundedness properties 
of solution set of Ric$(X)\le0$ for controllable and uncontrollable systems in the following section. 

 \section{Boundedness/unboundedness of the solution set of ARI}\label{ari}
%
 In this section, we study boundedness properties of the set of solutions of ARI and how they are related to the position of eigenvalues of the 
 feedback matrix for controllable and uncontrollable systems and the associated invariant subspaces. 
 We show that controllability is a necessary and sufficient condition for boundedness of the solution set of ARI (Theorem $\ref{bdns}$) under Assumption 
 $\ref{assum1}$. In the next two subsections, we observe the effect of the position of eigenvalues of the system matrix on the solution set of ARI 
 for controllable and uncontrollable systems. 

 \subsection{$(A_0,B)$ controllable}
   Suppose $(A_0,B)$ is controllable. Let $D_J$ be the upper triangular Jordan cannonical form associated with $A_0^T$. By $D_{J_{k}}$, we denote a $k\times k$ sub-matrix of $D_J$ which 
   is also in the Jordan form.  We now give the maximum and the minimum solution among all the solutions of the reduced order simplified $k \times k$
ARI: $-D_{J_k}\mathcal{L}-\mathcal{L}D_{J_k}^T+\mathcal{L}M_k\mathcal{L}\leq0$ where $(1\le k \le n)$. 
 First, we consider a special case where all eigenvalues of $A_0$ lie completely either in the left half complex plane or the right half complex plane 
 (Lemma \ref{maxsol}) followed  by the general case (Theorem \ref{maxminD}). 
 \begin{lemma}\label{pod}
  Suppose $(A_0,B)$ is controllable and $D_{J}$ be the Jordan canonical form associated with $A_0^T$ such that Spec$(D_{J})$ lies in
the open right half complex plane. Let the unique rank $n$ solution of the simplified ARE:
\begin{eqnarray}
 -D_{J}\mathcal{L}-\mathcal{L}D_{J}^T+\mathcal{L}M\mathcal{L}=0\nonumber
\end{eqnarray}
 be denoted by $\mathcal{L}^*$. If 
 $\bar{\mathcal{L}}$ is a rank $k$ solution $(k<n)$ of the simplified ARE, then $\bar{\mathcal{L}}\leq \mathcal{L}^*$.
 \end{lemma}
\begin{pf}
 Let $\bar{\mathcal{L}}$ be a rank $k$ solution of the simplified ARE: $-D_{J}\mathcal{L}-\mathcal{L}D_{J}^T+\mathcal{L}M\mathcal{L}=0$. 
 The column span of $\bar{\mathcal{L}}$ must be $D_J-$invariant. We can have an ordered basis of $\mathbb{R}^n$ whose first $k$ 
 vectors are defined by $k$ linearly independent columns of $\bar{\mathcal{L}}$. For the ease of notation, we continue to denote 
 both $\mathcal{L}^*$ and $\bar{\mathcal{L}}$ by the same notation after the change of basis.
 
 We know that (\cite{are} Theorem $13$, Corollary $13.1$ and other results of Section $4$), Schur complement of appropriate principal sub-matrices of $\mathcal{L}^*$ give low rank solutions of 
  simplified ARE: $-D_{J}\mathcal{L}-\mathcal{L}D_{J}^T+\mathcal{L}M\mathcal{L}=0$. 
  Let $\mathcal{L}^* = \left[\begin{array}{cc}\mathcal{L}^{11}&\mathcal{L}^{12}\\(\mathcal{L}^{12})^T&\mathcal{L}^{22}\end{array}\right]$ 
  (where $\mathcal{L}^{11}$ is $k\times k$ block and remaining blocks are of appropriate dimension). Since 
  all eigenvalues of $D_J$ are in the open right half plane, $\mathcal{L}^*>0$ (\cite{snyzak}). Therefore, $\mathcal{L}^{22} >0$. Let $\mathcal{L}^{22}=NN^T$ 
  where $N$ is an $(n-k)\times(n-k)$ invertible matrix. Therefore, 
 \begin{eqnarray}
  \bar{\mathcal{L}}&=& \left[\begin{array}{cc}\mathcal{L}^{11}-\mathcal{L}^{12}(\mathcal{L}^{22})^{-1}(\mathcal{L}^{12})^T&0\\0&0\end{array}\right].\nonumber\\
  \Rightarrow \mathcal{L}^* -\bar{\mathcal{L}} &=& \left[\begin{array}{cc}\mathcal{L}^{12}(\mathcal{L}^{22})^{-1}(\mathcal{L}^{12})^T&\mathcal{L}^{12}\\(\mathcal{L}^{12})^T&\mathcal{L}^{22}\end{array}\right]\label{eq_rk}\\
  & =&\left[\begin{array}{c}\mathcal{L}^{12}(N^{-1})^T\\N\end{array}\right]\left[\begin{array}{cc}N^{-1}(\mathcal{L}^{12})^T&N^T\end{array}\right] \geq0.\label{poseq}
  \end{eqnarray} \qed
\end{pf}

If in the above lemma, $D_J$ has all eigenvalues in the open left half plane, then we have $\bar{\mathcal{L}}\geq \mathcal{L}^*$. 
\begin{lemma}\label{lembdcom}
Suppose $-D_{J_k}$ has all eigenvalues in the open left half plane $(1\le k \le n)$. Then, $P= \int_{0}^{\infty}e^{-D_{J_k}^Tt}Ce^{-D_{J_k}t}dt$ is the unique 
solution of the Lyapunov equation 
 $P (-D_{J_k})+(-D_{J_k}^T)P = -C$.
\end{lemma}
\begin{pf}
 \begin{eqnarray}
 P (-D_{J_k})+(-D_{J_k}^T)P&= &\int_{0}^{\infty}e^{-D_{J_k}^Tt}Ce^{-D_{J_k}t}(-D_{J_k})dt+ \nonumber\\
 &&\int_{0}^{\infty}(-D_{J_k}^T)e^{-D_{J_k}^Tt}Ce^{-D_{J_k}t}dt\nonumber\\
 &=&\int_{0}^{\infty}\frac{d}{dt}e^{-D_{J_k}^Tt}Ce^{-D_{J_k}t}dt =-C.\nonumber
\end{eqnarray}
It can be shown that $P$ is the unique solution of this Lyapunov equation (since $D_{J_k}$ and $-D_{J_k}^T$ have no common eigenvalues (Theorem $4.4.6$ of \cite{hornj}).\qed
\end{pf}

\begin{lemma}\label{maxsol}
 Suppose $(A_0,B)$ is controllable and $D_{J}$ be the Jordan canonical form associated with $A_0^T$ such that Spec$(D_{J})$ lies in
the open right half complex plane.  
 Let $D_{J_k}$ be a $k \times k$ $(1\le k \le n)$ sub-matrix of $D_J$ in the Jordan canonical form. Let the unique rank $k$ solution of the reduced order 
 simplified ARE:
 \begin{eqnarray}
  -D_{J_k}\mathcal{L}-\mathcal{L}D_{J_k}^T+\mathcal{L}M_k\mathcal{L}=0\nonumber
 \end{eqnarray}
 be denoted by $\mathcal{L}_k^*$. Then every solution $\hat{\mathcal{L}}$ of the reduced order 
 simplified ARI: 
 \begin{eqnarray}
 -D_{J_k}\mathcal{L}-\mathcal{L}D_{J_k}^T+\mathcal{L}M_k\mathcal{L}\leq0\nonumber 
 \end{eqnarray}
 is such that $0\leq\hat{\mathcal{L}}\leq\mathcal{L}_k^*$.
\end{lemma}
\begin{pf}
 Suppose $\hat{\mathcal{L}}$ is a rank $k$ matrix that satisfies the reduced order simplified ARI:
 $-D_{J_k}\mathcal{L}-\mathcal{L}D_{J_k}^T+\mathcal{L}M_k\mathcal{L}\leq0$. 
 (By adding a positive semidefinite matrix $\mathcal{L}M'\mathcal{L}$ to the above ARI, we can convert it into an appropriate reduced order simplified ARE and applying 
 Theorem $12$ of \cite{are}, there exists a rank $k$ solution to the ARE which also gives a rank $k$ solution to the reduced order ARI considered above.) 
 First, we prove that $\hat{\mathcal{L}}\leq\mathcal{L}_k^*$. 
 Let $Y_k^*={\mathcal{L}_k^*}^{-1}$ and $\hat{Y}=\hat{\mathcal{L}}^{-1}$.
Let $\Delta=\hat{Y}-Y_k^*$. Pre and post-multiplying the ARI: $-D_{J_k}\hat{\mathcal{L}}-\hat{\mathcal{L}}D_{J_k}^T+\hat{\mathcal{L}}M_k\hat{\mathcal{L}}\le0$ with $\hat{Y}$, 
one obtains 
\begin{eqnarray}
 -\hat{Y}D_{J_k}-D_{J_k}^T\hat{Y}+M_k\leq0 \Rightarrow -{Y_k}^*D_{J_k}-D_{J_k}^T{Y_k}^*+M_k-\Delta D_{J_k}-D_{J_k}^T\Delta\leq0. \nonumber 
\end{eqnarray}
Now using $-{Y_k}^*D_{J_k}-D_{J_k}^T{Y_k}^*+M_k=0$, we obtain: 
\begin{eqnarray}
 -\Delta D_{J_k}-D_{J_k}^T\Delta \leq 0.\nonumber
\end{eqnarray}
 Let $-\Delta D_{J_k}-D_{J_k}^T\Delta = -C$ where $C\ge 0$. We can write this equation as
$\Delta (-D_{J_k})+(-D_{J_k}^T)\Delta = -C$. 
From Lemma $\ref{lembdcom}$, $\Delta=P=\int_{0}^{\infty}e^{-D_{J_k}^Tt}Ce^{-D_{J_k}t}dt$ is the unique solution of $\Delta (-D_{J_k})+(-D_{J_k}^T)\Delta = -C$. 
Observe that $P\ge0$ therefore, $\Delta \ge 0$. This implies that $\hat{Y}\ge Y_k^*\Rightarrow\hat{\mathcal{L}} \le \mathcal{L}_k^*$.

Now suppose that $\hat{\mathcal{L}}$ has rank strictly less than $k$ and it satisfies the reduced order simplified
 ARI: $-D_{J_k}\mathcal{L}-\mathcal{L}D_{J_k}^T+\mathcal{L}M_k\mathcal{L}\leq0$. From the results in the earlier sections, we know that the column span
 of $\hat{\mathcal{L}}$ must be a $D_{J_k}$-invariant subspace. Therefore, it is enough to consider the appropriate principal sub-matrix of the
 reduced order simplified ARI, as the $D_{J_k}$-invariant subspaces are spanned by a collection of elementary basis vectors $e_i$s with
 $i \in I \subset \{1,2, \cdots, k\}$. The above proof now shows that $\hat{\mathcal{L}} \le \mathcal{L}^*_S$, where $\mathcal{L}^*_S$ is the
 Schur complement of $\mathcal{L}_k^*$ with respect to the principal sub-matrix indexed by the set $\{1,2,\cdots,k\} \setminus I$. But then, 
 $\mathcal{L}^*_S \le \mathcal{L}_k^*$ (Lemma \ref{pod}) and therefore, $\hat{\mathcal{L}} \le \mathcal{L}_k^*$.

 Now we prove the other inequality. Again suppose $\hat{\mathcal{L}}$ is full rank and $\hat{Y}=\hat{\mathcal{L}}^{-1}$. 
 Again pre and post-multiplying the ARI: $-D_{J_k}\hat{\mathcal{L}}-\hat{\mathcal{L}}D_{J_k}^T+\hat{\mathcal{L}}M_k\hat{\mathcal{L}}\le0$ with $\hat{Y}$, we 
 obtain: 
 $-\hat{Y}D_{J_k}-D_{J_k}^T\hat{Y}+M_k\le0$. Hence, $-\hat{Y}D_{J_k}-D_{J_k}^T\hat{Y}\le -M_k$ where $M_k\ge 0$. Let $-\hat{Y}D_{J_k}-D_{J_k}^T\hat{Y} = -F$ where
 $F\ge0$. Again using Lemma $\ref{lembdcom}$, we have $\hat{Y}=\int_{0}^{\infty}e^{-D_{J_k}^Tt}Fe^{-D_{J_k}t}dt$. Thus, $\hat{Y}\ge0$
 (since $F\ge0$) and since it is full rank, $\hat{Y}^{-1}=\hat{\mathcal{L}}>0$.

 If  $\hat{\mathcal{L}}$ is  not full rank then restricting to the lower rank case and applying the similar arguments above, we obtain $\hat{\mathcal{L}}\ge0$.\qed
\end{pf}

Observe that if Spec$(D_{J_k})$ lies in the open left half plane, then following the same set of arguments as Lemma \ref{maxsol}, one can
conclude that $0 \ge \hat{\mathcal{L}} \ge \mathcal{L}^*$; where $\mathcal{L}^*$ is the unique rank $k$ solution of the corresponding ARE.
Note that the result above implies that solutions of Ric$(X)\le0$ satisfy $0\le X \le X^*$ when $A_0$ has all eigenvalues in the open
right half complex plane where $X^*$ is the full rank solution of ARE. Therefore, one can obtain \emph{Willems' result: 
$K_{min}\le K \le K_{max}$} (\cite{willems}, \cite{schererII}) for ARI by using $K=K_0+X$.

One can now combine these results along with the results of the earlier sections to obtain results about the most general case of the simplified ARI:
$-D_J\mathcal{L}-\mathcal{L}D_J^T+\mathcal{L}M\mathcal{L} \le 0$. 
Assume that $D_J$ has nonzero eigenvalues which are not purely imaginary. Without loss of generality, we 
assume that $D_J =$ diag$(D_{J_r}, D_{J_\ell})$ where $D_{J_r}\in \mathbb{R}^{k \times k}$ contains all the eigenvalues in the open right half plane and $D_{J_\ell}\in \mathbb{R}^{(n-k) \times (n-k)}$ contains all the eigenvalues in the
open left half plane. 
Let $\mathcal{L}^*$ be a maximal rank solution of the simplified ARE: $-D_J\mathcal{L}-\mathcal{L}D_J^T+\mathcal{L}M\mathcal{L} = 0$. 
Let $\mathcal{L}_r$ be the rank $k$ solution of the $k\times k$ reduced order simplified ARE: $-D_{J_r}\mathcal{L}-\mathcal{L}D_{J_r}^T+\mathcal{L}M_k\mathcal{L} = 0$ and 
 $\mathcal{L}_{\ell}$ be the rank $(n-k)$ solution of the $(n-k)\times (n-k)$ reduced order simplified ARE: $-D_{J_{\ell}}\mathcal{L}-\mathcal{L}D_{J_{\ell}}^T+\mathcal{L}M_{(n-k)}\mathcal{L} = 0$. 
From the results of \cite{are} (Section $4$), 
we know that $\mathcal{L}_r$ and $\mathcal{L}_{\ell}$ can be obtained from 
Schur complements of $\mathcal{L}^*$ with respect
to eigenvalues corresponding to the left/right half complex plane respectively. 
Note that 
$\mathcal{L}_r^* =$ diag$(\mathcal{L}_r, 0)$ and $\mathcal{L}_{\ell}^* =$ diag$(0,\mathcal{L}_{\ell})$ are solutions of the simplified 
ARE: $-D_J\mathcal{L}-\mathcal{L}D_J^T+\mathcal{L}M\mathcal{L} = 0$. This brings us to an important theorem as follows.

\begin{theorem}\label{maxminD}
Suppose $(A_0,B)$ is controllable and $A_0$ has no purely imaginary or zero eigenvalues. Consider $D_J$ = diag$(D_{J_r},D_{J_\ell})$ which is the 
Jordan form of $A_0^T$ where $D_{J_r}$ contains all the eigenvalues of $D_J$ in the open right half plane and $D_{J_\ell}$ contains all the eigenvalues of $D_J$ in the
open left half plane. 
Let $\mathcal{L}^*$ be a maximal rank solution of the simplified ARE and let $\mathcal{L}_r^*$ and
$\mathcal{L}_{\ell}^*$ be the solutions of the simplified ARE 
\begin{eqnarray}
 -D_J\mathcal{L}-\mathcal{L}D_J^T+\mathcal{L}M\mathcal{L} = 0\nonumber
\end{eqnarray}
obtained by Schur complements of principal sub-matrices of $\mathcal{L}^*$ with respect
to modes corresponding to the left/right half complex planes respectively. Then every solution
$\hat{\mathcal{L}}$ of the simplified ARI: 
\begin{eqnarray}
 -D_J\mathcal{L}-\mathcal{L}D_J^T+\mathcal{L}M\mathcal{L} \le 0\nonumber
\end{eqnarray}
 is such that $\mathcal{L}_{\ell}^* \le \hat{\mathcal{L}} \le \mathcal{L}_r^*$.
\end{theorem}

\begin{pf}
 Partition the solution $\hat{\mathcal{L}} = \left[ \begin{array}{cc} L_{11} & L_{12} \\ L_{12}^T & L_{22} \end{array} \right]$ and
 $M=\left[ \begin{array}{cc} M_{11} & M_{12} \\ M_{12}^T & M_{22} \end{array} \right]$, that matches
 the partition of $D_J=\left[ \begin{array}{cc} D_{J_r} & 0 \\ 0 & D_{J_\ell} \end{array} \right]$. Due to controllability, $M_{11}\neq0$ and $M_{22}\neq0$. 
 Let $D_{J_r},L_{11},M_{11}\in \mathbb{R}^{k \times k}$ and  
 $D_{J_\ell},L_{22},M_{22}\in \mathbb{R}^{(n-k) \times (n-k)}$. Let $\mathcal{L}^*
  = \left[ \begin{array}{cc} L_{11}^* & L_{12}^* \\ (L_{12}^*)^T & L_{22}^* \end{array} \right]$ be the maximal rank solution of the 
  simplified ARE: $-D_J\mathcal{L}-\mathcal{L}D_J^T+\mathcal{L}M\mathcal{L} = 0$. Let $\mathcal{L}_r$ 
  be the full rank solution of the reduced order ARE: $-D_{J_r}\mathcal{L}-\mathcal{L}D_{J_r}^T+\mathcal{L}M_{11}\mathcal{L}=0$ and 
   $\mathcal{L}_\ell$ be the full rank solution of the reduced order ARE: $-D_{J_\ell}\mathcal{L}-\mathcal{L}D_{J_\ell}^T+\mathcal{L}M_{22}\mathcal{L}=0$. 
 Clearly, $\mathcal{L}_r^* =$ diag$(\mathcal{L}_r, 0)$ satisfies the simplified ARE $-D_{J}\mathcal{L}-\mathcal{L}D_{J}^T+\mathcal{L}M\mathcal{L}=0$ 
 where $\mathcal{L}_r$ is obtained from the maximal rank solution $\mathcal{L}^*
  = \left[ \begin{array}{cc} L_{11}^* & L_{12}^* \\ (L_{12}^*)^T & L_{22}^* \end{array} \right]$ by taking the Schur complement  
 with respect to $L_{22}^*$ (\cite{are} Section $4$). Similarly, $\mathcal{L}_{\ell}^* =$ diag$(0,\mathcal{L}_{\ell})$ 
 satisfies the simplified ARE $-D_{J}\mathcal{L}-\mathcal{L}D_{J}^T+\mathcal{L}M\mathcal{L}=0$. 
   
 Suppose $\hat{\mathcal{L}}$ is a full rank matrix and satisfies the strict ARI.  
 We use the method of Schur complement in the proof, for which we 
 first need to prove that blocks $L_{11}$ and $L_{22}$ of $\hat{\mathcal{L}}$ are invertible. 
 Consider the $(1,1)$ block of the ARI $-D_J\mathcal{L}-\mathcal{L}D_J^T+\mathcal{L}M\mathcal{L} \le 0$. Thus, 
 \begin{eqnarray}
  -D_{J_r}L_{11}-L_{11}D_{J_r}^T \leq -Z\nonumber
 \end{eqnarray}
 where $Z = \left[ \begin{array}{cc} L_{11} & L_{12} \end{array} \right]M\left[ \begin{array}{c} L_{11} \\ L_{12}^T \end{array} \right]\ge0$.
 Therefore, by Lemma $\ref{lembdcom}$, $L_{11} \ge 0$. Similarly, from the $(2,2)$ block of the simplified ARI,
 one can conclude that $L_{22} \le 0$.
%

 We show that $L_{11}$ is invertible. Suppose it is not invertible and let $v\in$ ker$(L_{11})$ 
 where $v \in \mathbb{R}^k$. 
 Note that $D_J\hat{\mathcal{L}}+\hat{\mathcal{L}}D_J^T-\hat{\mathcal{L}}M\hat{\mathcal{L}} >0$. Restricting to $(1,1)$ block, 
 we get
 \begin{eqnarray}
  D_{J_r}L_{11}+L_{11}D_{J_r}^T-L_{11}M_{11}L_{11}-L_{11}M_{12}L_{12}^T-\nonumber\\
  L_{12}M_{12}^TL_{11}-L_{12}M_{22}L_{12}^T >0. \label{ricinq}
 \end{eqnarray}
Pre and post multiplying the Inequality (\ref{ricinq}) by $v^T$ and $v$, we get $-v^TL_{12}M_{22}L_{12}^Tv >0$ which is not possible since 
$M_{22}\geq0$. Therefore, $L_{11}$ must be a full rank matrix. Similarly, $L_{22}$ is a full rank matrix. Thus, both $L_{11}$ and $L_{22}$ are invertible. 

 Now we prove the inequality of the theorem when $\hat{\mathcal{L}}$ is a full rank matrix. Let $\hat{Y} = \hat{\mathcal{L}}^{-1}$ and $Y_{11}$ be its $(1,1)$ block. 
 Observe that $Y_{11}$ satisfies $-Y_{11}D_{J_r}-D_{J_r}^TY_{11}+M_{11}<0$.
 Note that using the Schur complement, $Y_{11}= (L_{11}-L_{12}L_{22}^{-1}L_{12}^T)^{-1}$ hence, $Y_{11}$ is invertible and $\mathcal{L}=Y_{11}^{-1}$
 satisfies $-D_{J_r}\mathcal{L}-\mathcal{L}D_{J_r}^T+\mathcal{L}M_{11}\mathcal{L}<0$.
 Therefore, due to Lemma \ref{maxsol},
 \begin{eqnarray}
  Y_{11}^{-1}= L_{11}-L_{12}L_{22}^{-1}L_{12}^T\le \mathcal{L}_r. \label{inim}
 \end{eqnarray}
 
 Now consider the matrix $\mathcal{L}_r^* - \hat{\mathcal{L}} = \left[ \begin{array}{cc} \mathcal{L}_r - L_{11} & - L_{12} \\ - L_{12}^T & - L_{22} \end{array} \right]$.
 As $L_{22} < 0$, taking the Schur complement of $\hat{\mathcal{L}}$ with respect to $L_{22}$, we get $0\le (\mathcal{L}_r-L_{11})+L_{12}L_{22}^{-1}L_{12}^T$ 
 (from Equation (\ref{inim})). Therefore, $\mathcal{L}_r^* - \hat{\mathcal{L}} \ge 0$. Similarly, one can argue that
 $\mathcal{L}_{\ell}^* - \hat{\mathcal{L}} \le 0$.

  Now suppose $\hat{\mathcal{L}}$ is a full rank matrix and satisfies the non-strict ARI. Let $(\mathcal{L}_m)$ be a sequence of rank $n$ matrices 
satisfying the strict ARI such that $\lim_{m \to \infty}\mathcal{L}_m=\hat{\mathcal{L}}$. 
Note that $\mathcal{L}_{\ell}^*\leq \mathcal{L}_m$ for all $m\geq1$. Taking limits on both 
sides as $m$ goes to infinity, we get $\mathcal{L}_{\ell}^*\leq \hat{\mathcal{L}}$. Similarly, $ \hat{\mathcal{L}} \le \mathcal{L}_r^*$. 
Thus, any $\hat{\mathcal{L}}$ of rank $n$ which is a solution of the simplified ARI satisfies the inequality of the theorem.
 
 Finally, we consider the case when $\hat{\mathcal{L}}$ is not full rank and satisfies the simplified ARI. 
 Hence, the column span of $\hat{\mathcal{L}}$ must be a $D_J-$invariant subspace.
 If $\hat{\mathcal{L}}$ has rank $k$, then it is enough to consider $k\times k$ sub-matrix $\hat{\mathcal{L}}_k$ of
 $\hat{\mathcal{L}}$ which is nonzero (all remaining entries of $\hat{\mathcal{L}}$ are zero because its columns form a $D_J-$invariant subspace).
 Thus, we have a reduced order simplified ARI for $k\times k$ case with the full rank solution $\hat{\mathcal{L}}_k$. Let $({\mathcal{L}}_r)_k$ and
 $({\mathcal{L}}_\ell)_k$ be maximum and minimum rank $k$ solutions. Therefore, using similar arguments used earlier,
 \begin{eqnarray}
  ({\mathcal{L}}_\ell)_k \leq \hat{\mathcal{L}}_k \leq({\mathcal{L}}_r)_k.\\
  \mbox{We also have } {\mathcal{L}}_\ell^*\leq({\mathcal{L}}_l)_k \mbox { and } ({\mathcal{L}}_r)_k \leq {\mathcal{L}}_r^*.\label{ine3}
 \end{eqnarray}
 (Inequality (\ref{ine3}) follows from Lemma \ref{maxsol}.) 
 The above two inequalities imply that ${\mathcal{L}}_\ell^* \leq \hat{\mathcal{L}}\leq{\mathcal{L}}_r^*$. \qed
\end{pf}
\begin{example}
Let 
 \begin{eqnarray}
  A =  \left[\begin{array}{ccc}1 & 0&0\\0 & 2&0\\0&0&-4\end{array}\right],
  B=\left[\begin{array}{c}1 \\1\\1\end{array}\right],Q=0.\nonumber
  \end{eqnarray}
  Since $A$ is diagonal, we can write $A=A_0=D_J$ and $M=BB^T$. Observe that $D_J$ has two eigenvalues in the RHP and one eigenvalue in the LHP. The full rank 
  solution of the simplified ARE $-D_J\mathcal{L}-\mathcal{L}D_J^T+\mathcal{L}M\mathcal{L} = 0$ is given by
  \begin{eqnarray}
 \mathcal{L}^*=  \left[\begin{array}{ccc}
    6.4800  & -4.8000 &   1.9200\\ 
   -4.8000 &   4.0000  & -3.2000\\
    1.9200  & -3.2000  & -0.3200\end{array}\right].\nonumber
  \end{eqnarray}

  By taking the Schur complement with respect to the lower $1\times 1$ principal sub-matrix of $\mathcal{L}^*$ associated with eigenvalue $-4$ in the LHP, 
  we get a rank two solution $\mathcal{L}^*_r$ as discussed in the above theorem. 
 \begin{eqnarray}
\mathcal{L}_r^*= \left[\begin{array}{ccc}18 & -24&0\\-24 & 36&0\\0&0&0\end{array}\right]. \nonumber
\end{eqnarray}
Now similarly, by taking the Schur complement with respect to the upper $2\times 2$ principal sub-matrix of $\mathcal{L}^*$ corresponding to eigenvalues 
in the RHP, we obtain a rank one solution given by
\begin{eqnarray}
\mathcal{L}_\ell^*= \left[\begin{array}{ccc}0 & 0&0\\0 & 0&0\\0&0&-8\end{array}\right].\nonumber
 \end{eqnarray}
 Let
 \begin{eqnarray}
  \mathcal{L}_1 = \left[\begin{array}{ccc}0.72 & 0&-1.92\\0 & 0&0\\-1.92&0&-2.88\end{array}\right].\nonumber
 \end{eqnarray}
 We observe that $\mathcal{L}_1$ satisfies
 $-D_J\mathcal{L}_1-\mathcal{L}_1D_J^T+\mathcal{L}_1M\mathcal{L}_1 = 0$. It turns out that $\mathcal{L}_\ell^*\leq \mathcal{L}_1\leq \mathcal{L}_r^*$. This is true for all
  solutions of the simplified ARE. 

 Consider  
  \begin{eqnarray}
   \hat{\mathcal{L}} = \left[\begin{array}{ccc}9.216 & -12 &  -0.5760\\-12&18&0\\-0.5760 &0&-2.4640\end{array}\right],\nonumber  
  \end{eqnarray}
  which satisfies the ARI $-D_J\hat{\mathcal{L}}-\hat{\mathcal{L}}D_J^T+\hat{\mathcal{L}}M\hat{\mathcal{L}}\leq 0$. We observe that 
$\hat{\mathcal{L}}$  also satisfies the inequality $\mathcal{L}_\ell^*\leq \hat{\mathcal{L}}\leq \mathcal{L}_r^*$.
\end{example}

When one goes back to the equation Ric$(X) \le 0$, the above theorem translates to the existence of a maximal and a minimal solution. 
Willems (\cite{willems}) proved that when $(A,B)$ is controllable,
all the solutions $K$ of ARI: $-A^TK-KA-Q+KBB^TK\leq0$ satisfy the inequality $K_{min}\leq K \leq K_{max}$. Assume that $(A,B)$ is controllable. Then, fixing 
some  solution $K_0$ of the ARE, one obtains $A_0 = A - BB^TK_0$. Consider the inequality Ric$(X) \le 0$ obtained using this data. Then, one reaches a situation where Theorem~\ref{maxminD} is applicable. 
Thus, $K_{min}$ and $K_{max}$ in Willems' result really comes from $\mathcal{L}_{\ell}^*$ and $\mathcal{L}_r^*$ respectively.

Observe that in Theorem~\ref{maxminD}, we assumed that the eigenvalues of $A_0$ do not lie on the imaginary axis. 
If we relax this condition, then one needs to consider a block structure of 
$D_J$ of the form diag$(D_{J_0}, D_{J_r},D_{J_\ell})$, where the sub-matrix $D_{J_0}$ contains all the purely imaginary and zero eigenvalues of $D_J$. 
Assuming controllability, this translates to the following: the sub-matrix of $M$ corresponding to the sub-matrix 
$D_{J_0}$ of $D_J$ is positive semidefinite. From results in Section~3 (specifically Theorem~\ref{zerocont}, Theorem~\ref{pureimag} and Theorem~\ref{pureimagjord}), one can therefore conclude that the 
corresponding block of $\mathcal{L}$ (a solution of the ARI) must be zero. Thus, Theorem  \ref{maxminD} holds when $A_0$ has purely imaginary and zero eigenvalues (which are controllable). 

This completes the description of the effect of the position of eigenvalues of $A_0$ on the solution set of ARI for controllable system. 
\subsection{$(A_0,B)$ uncontrollable}
Now we focus on uncontrollable systems and the behavior of the solution set of ARIs. 
 We assume that a given $(A_0,B)$ pair is uncontrollable and demonstrate two cases where the solution set of ARI becomes unbounded. 
\begin{theorem}\label{unc}
 Suppose $(A_0,B)$ is uncontrollable such that one of the following holds:
 \begin{itemize}
  \item $A_0$ has a real uncontrollable eigenvalue $\lambda$. 
  \item $A_0$ has a pair of complex conjugate eigenvalues which are uncontrollable. 
 \end{itemize}
 
 Then, the solution set of Ric$(X)\le0$ becomes unbounded.
\end{theorem}
\begin{pf}
 Suppose $A_0$ has a real uncontrollable eigenvalue $\lambda$. Therefore, there exists $v^T$ such that $v^T[\lambda I - A_0\; B]=0$. Thus, $v^TB=0$ 
 and $A_0^Tv= \lambda v$. Using $X=\alpha vv^T$, we obtain Ric$(X)=-2\alpha\lambda vv^T$. Thus, choosing $\alpha$ according to the sign of $\lambda$, we have 
 infinitely many solutions of Ric$(X)\le0$ and the solution set becomes unbounded. 
 
 The case of complex conjugate eigenvalues follows similarly using $X=L\mathcal{L}L^T$ (where the columns of $L$ form an invariant subspace associated 
 with the complex conjugate eigenvalues) and the fact that $L^TB=0$. \qed
\end{pf}

The above theorem along with Corollary $\ref{zerouncont}$ and Corollary $\ref{uncpuim}$ allows us to conclude that \emph{for uncontrollable systems, the 
solution set of ARIs become unbounded}. Thus, as a consequence, we have the following result which is one of the main results of this paper. 
\begin{theorem}\label{bdns}
 Suppose a solution to the ARE $-A^TK-KA-Q+KBB^TK=0$ exists. Then $(A,B)$ is controllable if and only if the set of solutions of 
 $-A^TK-KA-Q+KBB^TK\leq0$ is bounded.
\end{theorem}
\begin{pf}
If $(A,B)$ is controllable, then from Theorem $\ref{maxminD}$, the solution set of $-A^TK-KA-Q+KBB^TK\leq0$ is bounded. If $A$ has uncontrollable 
eigenvalues, then from Theorem \ref{unc}, Corollary \ref{zerouncont} and Corollary \ref{uncpuim}, the solution set of $-A^TK-KA-Q+KBB^TK\leq0$ becomes unbounded. 
 \qed
\end{pf}

It turns out that for uncontrollable systems, one can further characterize the unboundedness properties depending on the position of uncontrollable 
eigenvalues in the complex plane which is stated in Theorem \ref{bdbelow}. A variant of this result
has appeared in \cite{schererII}, \cite{palbel}. 
We now state a result from \cite{are} which is used to prove Theorem \ref{bdbelow}.
\begin{theorem}\label{comeig}
 Suppose $(A_0,B)$ is controllable. Let $X_*$ (a solution of Ric$(X) = 0$) be a rank $k$ matrix, such that the column span of $X_*$ is an $A_0^T-$invariant subspace 
 corresponding to $k$ real eigenvalues of $A_0$ , say $\lambda_1 ,\ldots, \lambda_k$. Then $A_0$ and $A_1=A_0-BB^TX_*$ have $(n - k)$ common
eigenvalues and $(n - k)$ common right eigenvectors. Furthermore, the column span of $X_*$ is a $A_1^T-$invariant subspace corresponding to 
eigenvalues $-\lambda_1,\ldots,-\lambda_k$ of $A_1$.
\end{theorem}
\begin{pf}
 Refer Theorem $14$ of \cite{are}. \qed
\end{pf}

This theorem also holds when $A_0$ has complex eigenvalues (Lemma $5$ of \cite{are}). 
\begin{remark}\label{constr}
 From the above theorem, it is clear that even if $A_0$ does 
not have all its eigenvalues in the open right half plane, one can move using an appropriate solution $X_*$ (of Ric$(X)=0$) to a solution $K_*=K_0+X$ of the original ARE 
$-A^TK-KA-Q+KBB^TK=0$ such that corresponding 
$A_*$ has all eigenvalues in the open right half complex plane. 
\end{remark}

\begin{theorem}\label{bdbelow}
 Let $(A_0,B)$ be uncontrollable. 
 \begin{itemize}
  \item If all uncontrollable eigenvalues lie in the open right half plane, then the solution set of Ric$(X)\leq0$ is bounded from below.
  \item If all uncontrollable eigenvalues
 lie in the open left half plane, then the solution set of Ric$(X)\leq0$ is bounded from above.
 \item If uncontrollable eigenvalues lie in both the half planes, then
 solution set of Ric$(X)\leq0$ is neither bounded below nor bounded above.
 \end{itemize}
   
\end{theorem}
\begin{pf}
 Using $X=L\mathcal{L}L^T$ where columns of $L$ are eigenvectors/generalized eigenvectors (or the real and imaginary parts of eigenvectors/generalized eigenvectors) 
 of $A_0^T$, Ric$(X)$ is reduced to $-D_J\mathcal{L}-\mathcal{L}D_J^T+\mathcal{L}M\mathcal{L}$ (Lemma \ref{smplfd}). 
 If all uncontrollable eigenvalues lie in the open right half plane, then we may assume that the real part of all controllable 
 eigenvalues of $D_J$ is also positive (Remark \ref{constr}). 
 Using arguments similar to those in the proof of Lemma~\ref{maxsol},  
 one can show that such solutions
 of Ric$(X)\leq0$ are positive semidefinite and hence bounded from below. Similarly, if all uncontrollable eigenvalues lie in the open left half plane,
 we can show that solution set of Ric$(X)\leq0$ are negative semidefinite and hence bounded from above. 
 
 If the uncontrollable eigenvalues lie in both half planes, then we 
 subdivide the uncontrollable eigenvalues  
 into those that lie in the open left half plane and those that lie in the open right half plane. Observe that solutions $X$ corresponding to 
 $A_0^T$-invariant subspaces associated with the uncontrollable eigenvalues in the open left/right half plane would be bounded from above/below and 
 therefore, the result follows. \qed
\end{pf}

Thus, we obtain a complete characterization of the solution set of ARI for both controllable and uncontrollable systems assuming that a solution of the 
corresponding ARE exists. Note that the sign controllability was not assumed in any of the stated results. 
\section{Parametrization of solutions of ARI}\label{para}
Now we give a parametrization of solutions of the simplified ARI when $(A_0,B)$ is controllable.  
\begin{theorem}\label{ARIchar}
Let $D_{J_k}$ be a $k \times k$ $(1\le k \le n)$ matrix in Jordan form such that Spec$(D_{J_k})$ lies in
the open right half complex plane. 
  Then, all rank $k$ solutions of the reduced order simplified strict ARI: 
  \begin{eqnarray}
   -D_{J_k}\mathcal{L}-\mathcal{L}D_{J_k}^T+\mathcal{L}M_k\mathcal{L} < 0\nonumber
  \end{eqnarray}
 are parametrized by all $k\times k$ positive definite matrices. 
 Furthermore, all rank $k$ solutions of the reduced order simplified ARI: 
  \begin{eqnarray}
  -D_{J_k}\mathcal{L}-\mathcal{L}D_{J_k}^T+\mathcal{L}M_k\mathcal{L} \leq 0\nonumber 
  \end{eqnarray}
 are parametrized by all $k\times k$ positive semidefinite matrices.
\end{theorem}
\begin{pf}
We show that any $k\times k$ positive definite matrix $P$ corresponds to a solution of the reduced order strict ARI and conversely. 
Let $P$ be a $k\times k$ positive definite matrix and consider the equation $-\Delta D_{J_k} - D_{J_k}^T\Delta = -P$. Since eigenvalues of $-D_{J_k}$ are
in the open left half plane, $\Delta$ is positive definite (Lemma $\ref{lembdcom}$).
Let the unique rank $k$ solution of the reduced order simplified ARE: $-D_k\mathcal{L}-\mathcal{L}D_k^T+\mathcal{L}M_k\mathcal{L}=0$ be denoted by $\mathcal{L}^*$
which is positive definite by Lemma \ref{maxsol}.
Let $Y^*=(\mathcal{L}^*)^{-1}$ and $\hat{Y}=Y^*+\Delta$. Clearly, $\hat{Y}$ is positive definite hence invertible.
Observe that $ -\hat{Y}D_{J_k}-D_{J_k}^T\hat{Y}+M_k = -\Delta D_{J_k} - D_{J_k}^T\Delta < 0$. Now pre and post multiplying
$ -\hat{Y}D_{J_k}-D_{J_k}^T\hat{Y}+M_k<0$ by $\hat{Y}^{-1} = \hat{\mathcal{L}}$, we get
$-D_{J_k}\hat{\mathcal{L}}-\hat{\mathcal{L}}D_{J_k}^T+\hat{\mathcal{L}}M_k\hat{\mathcal{L}} < 0$. Thus, all $k\times k$ positive definite matrices $P$
correspond to rank $k$ solutions satisfying the reduced order simplified ARI with strict inequality.

Now we want to show that associated with each solution of the reduced order strict ARI, we have a positive definite matrix $P$. 
Suppose $\hat{\mathcal{L}}$ is a rank $k$ solution satisfying the strict inequality $-D_{J_k}\hat{\mathcal{L}}-\hat{\mathcal{L}}D_{J_k}^T+
\hat{\mathcal{L}}M_k\hat{\mathcal{L}} < 0$. 
 By Lemma \ref{maxsol}, $\hat{\mathcal{L}} \le\mathcal{L}^*$.
Let $Y^*=(\mathcal{L}^*)^{-1}$ and $\hat{Y} = \hat{\mathcal{L}}^{-1}$.  Therefore, $\hat{Y}\geq Y^*$. Let 
$\hat{Y}=Y^*+\Delta$ where $\Delta \geq0$. Observe that $ -\hat{Y}D_{J_k}-D_{J_k}^T\hat{Y}+M_k <0$ which implies that
 $-\Delta D_{J_k} - D_{J_k}^T\Delta < 0$. 
Therefore, $\Delta$ satisfies $-\Delta D_{J_k} - D_{J_k}^T\Delta = -P$ where $P$ is a positive definite matrix. 
 Thus, all rank $k$ solutions satisfying strict simplified ARI are parametrized by $k\times k$ positive definite matrices.
 
The second statement follows using identical arguments above by using the non-strict inequality in place of the strict inequality and positive semidefinite 
matrices in place of positive definite matrices. 
\qed
\end{pf}

Note that the above theorem holds for all $1\le k \le n$ and the 
 low rank solutions can be parametrized by restricting the ARI to the lower rank case. 
Results similar to  Theorem~\ref{ARIchar} can be obtained for the case when Spec$(D_J)$ lies in the open left half complex plane. 
Parametrization of solutions of an ARI for any general case can be reduced to the above case by Remark \ref{constr}. 

A parametrization of solutions of an ARI has also appeared in 
\cite{schererII}, \cite{ferrIV}, \cite{pav}, \cite{fau} which is different from our approach of rank parametrization. Moreover, a 
parametrization of solutions of ARE has appeared in \cite{mehr3}, \cite{wim4}.

 \section{Conclusion}
  We gave a complete characterization of solutions of ARI: Ric$(X)\leq0$ and obtained a rank parametrization solutions
 of Ric$(X)\leq0$ (Theorem~\ref{ARIchar}) assuming that a solution of the corresponding ARE exists (Assumption~\ref{assum1}). 
 First, we translated the ARI problem 
 into an equivalent problem 
 involving a homogeneous inequality: Ric$(X)\leq0$ by fixing a solution of the corresponding ARE.  
 We proved that for controllable systems, if the feedback matrix $A_0$ has all eigenvalues in the open right half plane, then solutions $\mathcal{L}$ of the simplified ARI
 $-D_J\mathcal{L}-\mathcal{L}D_J^T+\mathcal{L}M_k\mathcal{L}\leq 0$ satisfy the inequality $0\le \mathcal{L} \le \mathcal{L}^*$ where
 $\mathcal{L}^*$ is the maximal solution. Similar result holds with inequalities reversed when eigenvalues of $A_0$ are in the open left half plane.
 We further showed that under certain conditions, a maximal rank solution of Ric$(X) = 0$ may be used to obtain an upper and lower bound 
 ($\mathcal{L}_{r}^*$ and $\mathcal{L}_{\ell}^*$ respectively) for  
 solutions of Ric$(X) \leq 0$ (Theorem~\ref{maxminD}). This in turn provides a classical result of Willems ($K_{min}\leq K \leq K_{max}$ where $K$ is a 
 solution of the ARI: $-A^TK-KA-Q+KBB^TK\leq0$) as a special case. 
Without assuming sign controllability, 
 we showed that the solution 
 set of the ARI is bounded if and only if the system is controllable (Theorem $\ref{bdns}$). 
 
 We further showed that if all uncontrollable eigenvalues are in the right half plane then the solution set of Ric$(X)\leq0$ is
 bounded from below and if they are in the left half plane, then it is bounded from above. If the uncontrollable eigenvalues lie in both half planes, then the
 solution set of Ric$(X)\leq0$ is neither bounded above nor bounded below. This observation along with the results for controllable systems 
 captures the relationship between the position of eigenvalues of $A_0$ and the solution set of ARI for controllable as well as uncontrollable systems.

\bibliographystyle{model2-names}
\bibliography{sanandb}







\end{document}